\documentclass[a4paper,12pt]{article}

\usepackage{amsmath,amsfonts}
\usepackage{mathrsfs} 
\usepackage{parskip}
\usepackage{enumerate}
\usepackage[T1]{fontenc}
\usepackage[utf8]{inputenc}
\usepackage{authblk}

\newtheorem{theorem}{Theorem}
\newtheorem{lemma}{Lemma}

\def\Pb{\mathbf{P}}

\def\Ex{\mathbf{E}}

\def\II{\mathbb{I}}
\def\JJ{\mathbb{J}}

\def\1{\mbox{1\hspace{-.25em}I}}
\begin{document}
\title{On  Parameter Estimation of Hidden Telegraph Process}
\author[1]{R. Z. Khasminskii}
\author[2]{Yu.A. Kutoyants}
\affil[1]{\small Wayne State University, Detroit, USA}
\affil[1]{\small Institute for Information Transmission Problems, Moscow, Russia}
\affil[2]{\small National Research University ``MPEI'', Moscow, Russia, }
\affil[2]{\small University of  Maine,  Le Mans,  France}

\date{}
\maketitle

\begin{abstract}
 The problem of parameter estimation by the observations of the
two-state telegraph process in the presence of white Gaussian noise is
considered. The properties of estimator of the method of moments are described
in the asymptotics of large samples. Then this estimator is used as
preliminary one to construct the one-step MLE-process, which provides the
asymptotically normal and asymptotically efficient estimation of the unknown
parameters.  
\end{abstract}


\section{Introduction}

This work is devoted to the problem of parameter estimation by the
observations in continuous time $X^T=\left(X\left(t\right),0\leq t\leq
T\right)$ of the following stochastic process
\begin{equation}
\label{01}
{\rm d}X_t=Y\left(t\right)\,{\rm d}t+{\rm d}W_t,\quad
X_0, 
\end{equation}
here $W_t,0\leq t\leq T$ is a standard Wiener process, $X\left(0\right)=X_0$
is the initial value independent of the Wiener process and $Y\left(t\right),0\leq
t\leq T$ is a two-state ($y_1$ and $y_2$) stationary Markov process with
 transition rate  matrix
\begin{align*}
\begin{pmatrix}
-\lambda  &\lambda \\
\mu &-\mu 
\end{pmatrix}
.
\end{align*}
We suppose that the values $\lambda>0 $ and $\mu>0 $ are unknown and we have to
estimate the two-dimensional parameter $\vartheta =\left(\lambda ,\mu
\right)\in \Theta $, here $\Theta =\left(c_0,c_1\right)\times\left(c_0,c_1\right) $
by the observations $X^t,0<t\leq T$, i.e., the estimator $\vartheta
_{t,T}^\star , 0<t\leq T$ is stochastic process. Here $c_0<c_1$ are positive constants. 

Therefore, our goal is to construct an {\it on-line} estimator-process
$\vartheta _{T}^\star =\left(\vartheta _{t,T}^\star, 0<t\leq T\right)$, which
can be sufficiently easy to evaluate and is asymptotically optimal in some
sense, as $T\rightarrow \infty $. This estimator-process we constructe in two
steps. First we introduce a learning interval $\left[0,T^\delta \right]$,
here $\delta \in (\frac{1}{2}, 1)$ and propose a $T^{\frac{\delta }{2}}
$-consistent preliminary estimator constructed of the method of moments.  Then
we improve it up to asymptotically efficient one with the help of slightly
modified one-step MLE procedure.

Such model of observations is called ``Hidden Markov Model'' (HMM) or
partially observed system. There exists an extensive study of such type HMM
for discrete time models of observations, see, for example, \cite{EAM95},
\cite{BRR98}, \cite{CMT05} and the references therein. For continuous time
observation models this problem is not so-well studied.  See, for example,
Elliot et al. \cite{EAM95}, here the Part III ``Continuous time estimation''
is devoted to such models of observations.  One can find there the discussion
of the problems of filtration and parameter estimation in the case of
observations of finite-state Markov process in the presence of White
noise. The problem most close to our statement was studied by Chigansky
\cite{PCh09}, who considered the parameter estimation in the case of hidden
finite-state Markov process by continuous time observations. He showed the
consistency, asymptotic normality and asymptotic efficiency of the MLE of
one-dimensional parameter. The case of two-state hidden telegraph process
studied in our work was presented there as example but there supposed that
$\lambda =\mu $.  The problem of parameter estimation for similar models
including the identification of partially observed linear processes (models of
Kalman filtration) were studied by many authors. Let us mention here
\cite{Kut84}, \cite{DZ86}, \cite{ZD88}, \cite{Kut94}. The problem of
asymptotically efficient estimation for the model of telegraph process
observed in discrete times was studied in \cite{IY09}. 

The proposed here one-step MLE-process is motivated by the work of Kamatani
and Uchida \cite{KU15} who introduced Newton-Raphson multi-step estimators in the problem of
parameter estimation by discrete time observations of diffusion process. In
particular, it was shown that multi-step Newton-Raphson procedure allows to
improve the rate of convergence of preliminary estimator up to asymptotically
efficient one. Note that preliminary estimator there is constructted by all
observations and they studied  estimator of the unknown parameter.
Applied in this work estimator-process uses the preliminary estimator
constructed by the observations on the initial learning interval and follows
the similar construction as that one  introduced in the work \cite{Kut15} (see as well
\cite{Kut14}).

\section{Problem statement and auxiliary results}

We start with description of the MLE for this model of observations. The
stochastic process \eqref{01} according to {\it innovation theorem} (see
\cite{LS01}, Theorem 7.12) admits the representation
\begin{equation*}
{\rm d}X_t=m\left(t,\vartheta \right)\,{\rm d}t+{\rm d}\bar
W_t, \quad X_0,\quad 0\leq t\leq T,
\end{equation*}
here $m\left(t,\vartheta \right) $ is the conditional expectation
\begin{equation*}
m\left(t,\vartheta \right)=\Ex_\vartheta \left[Y\left(t\right)|{\cal
    F}_t^X\right]= y_1\Pb_\vartheta \left(Y\left(t\right)=y_1|{\cal F}_t^X
\right)+y_2\Pb_\vartheta \left(Y\left(t\right)=y_2|{\cal F}_t^X \right) .
\end{equation*}
Here ${\cal F}_t^X $ is the  $\sigma $-algebra generated by the  observations up
to time $t$, i.e., ${\cal F}_t^X :=\sigma \left( X_t,0\leq s\leq
t\right) $ and $\bar W_t,0\leq t\leq T$ is an {\it innovation}
Wiener process. 

Let us denote 
\begin{equation*}
\pi \left(t,\vartheta \right)=\Pb_\vartheta \left(Y\left(t\right)=y_1|{\cal F}_t^X
\right) ,\qquad \Pb_\vartheta \left(Y\left(t\right)=y_2|{\cal F}_t^X
\right)=1-\pi \left(t,\vartheta \right). 
\end{equation*}
Hence
\begin{align*}
m\left(t,\vartheta \right)=y_2+\left(y_1-y_2\right)\pi \left(t,\vartheta \right).
\end{align*}
The random process $\pi \left(t,\vartheta \right),0\leq t\leq T $ satisfies
the following equation (see \cite{LS01}, Theorem 9.1 and equation (9.23) there)
\begin{align}
{\rm d}\pi \left(t,\vartheta \right)&=\left[\mu -\left(\lambda +\mu \right)\pi
  \left(t,\vartheta \right)\right.\nonumber\\
&\quad \left. +\pi \left(t,\vartheta \right)\left(1-\pi
\left(t,\vartheta \right)\right)
\left(y_2-y_1\right)\left(y_2+\left(y_1-y_2\right)\pi \left(t,\vartheta
\right) \right)\right]{\rm d}t\nonumber\\
&\quad +\pi \left(t,\vartheta \right)\left(1-\pi
\left(t,\vartheta \right)\right)\left(y_1-y_2\right)\;{\rm d}X_t.
\label{03}
\end{align}
Denote by  $\left\{\Pb_\vartheta ^{\left(t\right)},\vartheta \in\Theta
\right\}$ the  measures induced by the observations $X^t=\left(X_s,0\leq s\leq t\right)$ of
stochastic processes \eqref{01} with different 
$\vartheta $ in the space of realizations ${\cal C}\left[0,t\right]$
(continuous on $\left[0,t\right]$ functions). These measures  are equivalent and the likelihood
ratio function 
\begin{align*}
L\left(\vartheta ,X^t\right)=\frac{{\rm d}\Pb_\vartheta ^{\left(t\right)}
}{{\rm d}\Pb_0 ^{\left(t\right)}}\left(X^t\right) ,\qquad \vartheta \in \Theta ,\quad 0<t\leq T
\end{align*}
can be written as follows
\begin{align*}
L\left(\vartheta ,X^t\right)=\exp\left\{\int_{0}^{t}m\left(s,\vartheta
\right){\rm d}X_s-\frac{1}{2}\int_{0}^{t}m\left(s,\vartheta
\right)^2{\rm d} s \right\}.
\end{align*}
Here $\Pb_0 ^{\left(t\right)} $ is the measure corresponding to $X^t$ with
$Y\left(s\right)\equiv 0$. 
 
The MLE-process $\hat\vartheta_{t,T}$ is defined by the equation 
\begin{align}
\label{5a}
L\left(\hat\vartheta_{t,T} ,X^t\right)=\sup_{\vartheta \in \Theta
}L\left(\vartheta ,X^t\right), \qquad 0<t\leq T.
\end{align}

It is known that in the one-dimensional case ($d=1, \lambda =\mu=\vartheta  $) the MLE
$\hat\vartheta_{T,T}=\hat\vartheta_{T} $ is
consistent, asymptotically normal 
$$
\sqrt{T}\left(\hat\vartheta_T-\vartheta  \right)\Longrightarrow {\cal
  N}\left(0,{\rm I}\left(\vartheta \right)^{-1}\right)
$$
and asymptotically efficient (see \cite{PCh09}, \cite{IH81}). Here ${\rm I}\left(\vartheta
\right) $ is the Fisher information. 

Note that the construction of the MLE-process $\hat\vartheta_{t,T},0<t\leq T$
according to \eqref{5a} and \eqref{03} is 
computationally difficult problem because we need to solve the family of
equations \eqref{03} for all $\vartheta \in \Theta $ and \eqref{5a} for all $t\in (0,T]$. 

We propose the following construction. First we study an estimator
$\bar\vartheta _T$ of the method of moments and show that this estimator is
$\sqrt{T}$-consistent, i.e.,
\begin{align*}
\Ex_\vartheta \left|\sqrt{T}\left(\bar\vartheta _T-\vartheta  \right)\right|^2\leq C.
\end{align*}
 Then using this estimator  $\bar\vartheta _{T^\delta }$  obtained by the
 observations on the {\it learning 
 interval} $\left[0,T^\delta \right]$, here $\frac{1}{2}<\delta <1$, we
 introduce the one-step MLE-process
\begin{align*}
\vartheta _{t,T}^\star=\bar\vartheta _{T^\delta }+T^{-1/2}{\II}_T\left(\bar\vartheta
_{T^\delta } \right)^{-1} \Delta _t\left(\bar\vartheta _{T^\delta },X^t\right),\qquad T^\delta \leq t\leq T.
\end{align*}
Here the empirical Fisher information matrix
\begin{align*}
{\II}_t\left(\vartheta
 \right)=\frac{1}{t}\int_{T^\delta }^{t}\dot m\left(\vartheta ,s\right)\dot
 m\left(\vartheta ,s\right)^* {\rm d}s\longrightarrow {\II}\left(\vartheta
 \right) ,
\end{align*}
as $t\rightarrow \infty , T^\delta =o\left(t\right)$ and the vector score-function process is
\begin{align*}
\Delta _t\left(\vartheta ,X^t\right)=\frac{1}{\sqrt{t}}\int_{T^\delta
}^{t}\dot m\left(\vartheta ,s\right) \,\left[{\rm d}X_s- m\left(\vartheta ,s\right){\rm d}s\right].
\end{align*}
Here and in the sequel dot means the derivation w.r.t. parameter $\vartheta $,
$\dot m\left(\vartheta ,t\right)$ is the vector-column of the derivatives
$\dot m_\lambda \left(\vartheta ,t\right) $ and $\dot m_\mu \left(\vartheta
,t\right) $. 
The estimator $\vartheta _{t,T}^\star $ is in some sense asymptotically
efficient. In particular for $\vartheta _{T,T}^\star=\vartheta _{T}^\star
$ we have
\begin{align*}
\sqrt{T}\left( \vartheta _{T}^\star-\vartheta \right)\Longrightarrow {\cal
  N}\left(0,{\II}\left(\vartheta \right)^{-1} \right),
\end{align*}
i.e., it is asymptotically equivalent to the MLE. Note that the calculation of
the estimator $\vartheta _{t,T}^\star $ for all $t\in \left[T^\delta
  ,T\right]$ requires the solution of the equation \eqref{03} for one value
$\vartheta =\bar \vartheta_{T^\delta } $ only. 

Recall as well the well-known properties of the Telegraph (stationary) process
$Y\left(t\right), t\geq 0$.
\begin{enumerate}
\item The stationary distribution of the process $Y\left(t\right)$ is
\begin{equation}
\label{1.2}
\Pb_\vartheta\left\{Y\left(t\right)=y_1\right\}=\frac{\mu }{\lambda +\mu },\qquad
\Pb_\vartheta\left\{Y\left(t\right)=y_2\right\}=\frac{\lambda  }{\lambda +\mu } 
\end{equation}
\item Let us denote
  $P_{ij}\left(t\right)=\Pb_\vartheta\left\{Y\left(t\right)=y_j|Y\left(0\right)=y_i\right\}$,\quad
  then solving the Kolmogorov equation we obtain
\begin{align}
P_{11}\left(t\right)&=\frac{\mu }{\lambda +\mu } +\frac{\lambda }{\lambda +\mu
}e^{-\left(\lambda +\mu \right)t}, \quad P_{12}\left(t\right)=\frac{\lambda  }{\lambda +\mu
} -\frac{\lambda }{\lambda +\mu }e^{-\left(\lambda +\mu \right)t} ,\nonumber\\
P_{21}\left(t\right)&=\frac{\mu }{\lambda +\mu } -\frac{\mu  }{\lambda +\mu
}e^{-\left(\lambda +\mu \right)t}, \quad P_{22}\left(t\right)=\frac{\lambda  }{\lambda +\mu
}+\frac{\mu }{\lambda +\mu }e^{-\left(\lambda +\mu \right)t} 
\label{1.3}
\end{align}
It follows from \eqref{1.2} and \eqref{1.3}  that
\begin{align}
K\left(s\right)&=\Ex_\vartheta Y\left(t\right)Y\left(t+s\right)=\left(\frac{y_1\mu
  +y_2\lambda }{\lambda +\mu }\right)^2 \nonumber\\
&+\left(y_2-y_1\right)^2\frac{\lambda
  \mu }{\left(\lambda +\mu \right)^2} e^{-\left(\lambda +\mu \right)s}
=\left(\bar Y\right)^2+ De^{-\left(\lambda +\mu \right)s},
\label{1.4}
\end{align}
here 
\begin{equation}
\label{def}
\bar Y=\Ex_\vartheta Y\left(t\right)=\frac{y_1\mu
  +y_2\lambda }{\lambda +\mu },\quad D=\left(y_2-y_1\right)^2\frac{\lambda
  \mu }{\left(\lambda +\mu \right)^2}.
\end{equation}

\item Let ${\cal F}_t^Y\subset {\cal F}$ be a family of $\sigma $-algebras,
  induced by the events $$\left\{ Y\left(s\right)=y_i,0\leq s\leq t,
  i=1,2\right\}.$$
Then  it follows from \eqref{1.3} that for some constant $K>0$ and $A<T$ and for
all $s>A,t>0$ the inequality 
\begin{equation}
\label{1.5}
\left|\Ex_\vartheta \left\{Y\left(s+T\right)Y\left(t+T\right)|{\cal F}_A^Y \right\} -\Ex_\vartheta
\left[Y\left(s\right)Y\left(t\right)\right]\right|<Ke^{-\left(\lambda +\mu \right)\left(T-A\right)}
\end{equation}
holds.
\end{enumerate}

\section{  Method of moments estimator}

Let us consider the problem of the construction of  $\sqrt{T}$-consistent estimators
of the parameter $\vartheta $ by the method of moments. Recall that we observe in
continuous 
time the stochastic process
\begin{equation}
\label{1.1}
{\rm d}X_t=Y\left(t\right){\rm d}t+{\rm d}W_t,\quad X_0,\quad  0\leq t\leq T,
\end{equation}
here $W_t, 0\leq t\leq T$ is a standard Wiener process, $X_0$ is independent
of $W_t,0\leq t\leq T$ initial value, 
$Y\left(t\right)=Y\left(t,\omega \right)$ is stationary Markov process 
with two states $y_1$ and $y_2$ and infinitesimal rate matrix
\begin{align*}
\begin{pmatrix}
-\lambda  &\lambda \\
\mu &-\mu 
\end{pmatrix}.
\end{align*}
The processes $Y\left(t \right),t\geq 0$ and $W_t,t\geq 0$ are independent.

We suppose for  simplicity that $T$ is an integer number. Introduce the  condition
\begin{equation}
\label{1.6}
\lambda \in
\left[c_0,c_1\right],\mu \in \left[c_0,c_1\right]
\end{equation}
here $c_0$ and $c_1$ are  some positive constants. 

To introduce  the estimators we need the following notations. 
\begin{itemize}
\item The function 
\begin{equation}
\label{fi}
\Phi \left(x\right)=\frac{1}{x}-\frac{1}{x^2}\left(1-e^{-x}\right).
\end{equation}
\item The statistics
\begin{align}
\label{2.2}
\zeta _T=\frac{1}{T}\sum_{i=0}^{T-1}\left[X_{i+1}-X_i\right]^2-1.
\end{align}
\item The random variable $\alpha _T$ is defined as a solution of the equation
\begin{align}
\label{2.11}
\zeta _T=\left(\frac{X_T}{T}\right)^2+2\eta _T\,\Phi \left(\alpha _T\right),
\end{align}
here 
\begin{align}
\label{2.12a}
\eta _T=\left(\frac{X_T}{T}-y_1\right)\left(y_2-\frac{X_T}{T}\right).
\end{align}
\item The event ${\cal A}_T$ that the equation \eqref{2.11} has a solution.
\item The random variable 
\begin{align}
\label{2.12}
\beta _T=\alpha _T\,\1_{\left\{{\cal
    A}_T\right\}}+\left(c_0+c_1\right)\,\1_{\left\{{\cal A}_T^c\right\}} 
\end{align}
\end{itemize}

The method of moments estimator is $\hat\vartheta _T=\left(\hat\lambda
_T,\hat\mu  _T\right)$,  here
\begin{equation}
\label{2.18}
\hat\lambda _T=\frac{\frac{X\left(T\right)}{T}-y_1}{y_2-y_1}\beta _T;\qquad
\hat\mu  _T=\frac{y_2-\frac{X\left(T\right)}{T}}{y_2-y_1}\beta _T .
\end{equation}

The properties of these estimators are given in the following theorem.

\begin{theorem}
\label{T1} Let the condition \eqref{1.6} holds. Then for the estimators
\eqref{2.18} and some constante $C>0$ we have for all $T>0$
\begin{equation}
\label{2.19}
\Ex_\vartheta \left[\sqrt{T}\left(\hat\lambda _T-\lambda \right)\right]^2<C,\qquad \Ex_\vartheta
\left[\sqrt{T}\left(\hat\mu _T-\mu  \right)\right]^2<C.
\end{equation}
\end{theorem}
The proof is given in several steps.

The next lemma gives $\sqrt{T}$-consistent estimator for $\bar Y$ (see \eqref{def}).
\begin{lemma}
\label{L1}
Let the condition \eqref{1.6} be fulfilled. Then  the estimator
$X_T/T$ is uniformly consistent for $\bar Y$ and for any $T>0$
\begin{align}
\label{con}
\Ex_\vartheta \left(\frac{X_T}{T}-\bar Y\right)^2\leq \frac{C}{T},
\end{align}
here the  constant $C>0$ does not depend on $\vartheta $.
\end{lemma}
{\bf Proof.} Making use \eqref{1.4} we obtain the fllowing relations
\begin{align*}
\Ex_\vartheta \left(\frac{X_T}{T}-\bar Y\right)^2&=\Ex_\vartheta\left|
\frac{1}{T}\int_{0}^{T}\left[Y\left(t\right)-\bar Y\right]{\rm d}t
+\frac{W_T}{T}\right|^2 \\
&=\frac{1}{T}+\frac{1}{T^2}\Ex_\vartheta\left|\int_{0}^{T}\left[Y\left(t\right)-\bar
  Y\right]{\rm d}t\right|^2\\ 
&\leq  \frac{1}{T}\left(   1+\frac{2\lambda \mu }{\left(\lambda +\mu\right)^3}
\left(y_2-y_1\right)^2 \right)\leq  \frac{1}{T}\left(   1+\frac{c_1^2 }{4c_0^3}
\left(y_2-y_1\right)^2 \right).
\end{align*}
{\it Corollary.} The existence of the
consistent estimator for $\frac{\lambda }{\lambda +\mu }$ and $\frac{\mu
}{\lambda +\mu }$  follows  from \eqref{1.2} and Lemma 1. Indeed, from the equality
\begin{align*}
\bar Y=\frac{\lambda }{\lambda +\mu }y_2+\frac{\mu
}{\lambda +\mu }y_1
\end{align*}
and Lemma 1 we obtain 
\begin{align}
&\Ex_\vartheta \left[\sqrt{T} \left(   \frac{T^{-1}X_T-y_1}{y_2-y_1}
  -\frac{\lambda  }{\lambda +\mu } \right)\right]^2 <C,\nonumber\\
&\Ex_\vartheta \left[\sqrt{T} \left(   \frac{y_2-T^{-1}X_T}{y_2-y_1}
  -\frac{\mu }{\lambda +\mu } \right)\right]^2 <C.
\label{1.7}
\end{align}
The statistics
\begin{align*}
\frac{X_T}{T}=\frac{1}{T}\int_{0}^{T}Y\left(t\right){\rm
  d}t+\frac{W_T}{T} 
\end{align*}
is the sum of a bounded a.s. random variable and an independent of it gaussian
random variable with parameters $\left(0,T^{-1}\right)$. Hence for $\eta _T$
defined in 
\eqref{2.12a} we can write the estimate
\begin{equation}
\label{1.8}
\Ex_\vartheta \left[\sqrt{T}\left(\eta  _T-D\right)\right]^2<C,
\end{equation}
here the constant $C>0$ does not depend on $T$ and $\vartheta $. The constant
$D$ is defined in \eqref{def}.

Note that from the condition \eqref{1.6} we have
\begin{align*}
\frac{\lambda \mu }{\left(\lambda +\mu \right)^2}>\frac{c_0^2}{4c_1^2}
\end{align*}
and  we easily obtain the estimate \eqref{1.8} for the estimator
\begin{equation}
\label{1.9}
\tilde \eta  _T=\max\left\{\eta  _T,\frac{c_0^2}{8c_1^2}\right\}
\end{equation}

\begin{lemma}
\label{L2} The following equality holds
\begin{equation}
\label{2.3}
\Ex_\vartheta \zeta _T=\bar Y^2+2D\Phi \left(\lambda +\mu \right)
\end{equation}
and under the condition \eqref{1.6} we have as well
\begin{equation}
\label{2.4}
\Ex_\vartheta\left[\sqrt{T}\left(\zeta _T-\Ex_\vartheta \zeta _T\right)\right]^2<C.
\end{equation}

\end{lemma}
{\bf Proof.} From stationarity of the process $Y\left(t\right)$ and \eqref{1.4}
we obtain
\begin{align}
\label{2.5}
\Ex_\vartheta \zeta _T&=\Ex_\vartheta \left[X_1-X_0\right]^2-1=\Ex_\vartheta
\int_{0}^{1} \int_{0}^{1}Y\left(s\right) Y\left(t\right){\rm d}s{\rm
  d}t+1-1\nonumber\\
&=\bar Y^2+2D\Phi \left(\lambda +\mu \right).
\end{align}
Denote
\begin{align*}
\gamma_i=\int_{i}^{i+1}Y\left(t\right){\rm d}t;\qquad \Delta
W\left(i\right)=W_{i+1}-W_i .
\end{align*}
Further, from the equality
\begin{align*}
\zeta _T-\Ex_\vartheta\zeta _T&=\frac{1}{T}\sum_{i=0}^{T-1}\left(\gamma_i^2-\Ex_\vartheta \gamma
_i^2\right) +\frac{2}{T}\sum_{i=0}^{T-1}\gamma_i\Delta W\left(i\right)
+\frac{1}{T}\sum_{i=0}^{T-1}\left(\Delta W\left(i\right)^2-1\right) 
\end{align*}
follows the estimate
\begin{align}
\label{2.6}
\Ex_\vartheta \left(\zeta _T-\Ex_\vartheta \zeta _T\right)^2&\leq \frac{3}{T^2}\Ex_\vartheta
\left(\sum_{i=0}^{T-1}\left(\gamma_i^2-\Ex_\vartheta \gamma 
_i^2\right)\right)^2 +\frac{12}{T^2}\Ex_\vartheta \left(\sum_{i=0}^{T-1}\gamma_i\Delta
W\left(i\right)\right)^2 \nonumber\\
& +\frac{3}{T^2}\Ex_\vartheta \left(  \sum_{i=0}^{T-1}\left(\Delta W\left(i\right)^2-1\right) \right)^2:=3J_1+12J_2+3J_3.
\end{align}
From stationarity of $Y\left(t\right)$ we obtain
\begin{align}
\label{2.7}
J_1&=\frac{1}{T^2}\Ex_\vartheta
\left(\sum_{i=0}^{T-1}\left(\gamma_i^2-\Ex_\vartheta \gamma
_i^2\right)\right)^2=\frac{1}{T^2}  \sum_{i=0}^{T-1}\sum_{j=0}^{T-1}  \Ex_\vartheta\left(\gamma_i^2-\Ex_\vartheta \gamma 
_i^2\right)\left(\gamma_j^2-\Ex_\vartheta \gamma
_j^2\right)\nonumber\\
&=\frac{1}{T^2}\sum_{i,j=0}^{T-1}\left\{\int_{0}^{1}\int_{0}^{1}\int_{0}^{1}\int_{0}^{1}\Ex_\vartheta \left\{
  Y\left(s\right)Y\left(t\right)\right.\nonumber\right.\\
&
  \qquad \left. \Ex_\vartheta\left[Y\left(\left|i-j\right|+s_1\right)Y\left(\left|i-j\right|+t_1\right)|{\cal
    F}_1^Y \right]\right\}{\rm d}s{\rm d}t{\rm d}s_1{\rm d}t_1-\Ex_\vartheta\gamma_0^2 \Bigr\}.
\end{align}
The estimate \eqref{1.5} allows to write
\begin{align*}
\left|\Ex_\vartheta\left[Y\left(\left|i-j\right|+s_1\right)Y\left(\left|i-j\right|+t_1\right)|{\cal
    F}_1^Y \right]-\Ex_\vartheta Y\left(s_1\right)Y\left(t_1\right)\right|\leq
K\,e^{-\left(\lambda +\mu \right)\left|j-i\right|} .
\end{align*}
From this estimate, \eqref{2.7} and \eqref{1.6} we obtain
\begin{align*}
J_1\leq \frac{K}{T^2}\sum_{i,j=0}^{T-1}e^{-\left(\lambda +\mu \right)\left|j-i\right|}\leq \frac{K_1}{T}.
\end{align*}
The following estimates are evident
\begin{align*}
J_2&=\frac{1}{T^2}\Ex_\vartheta \left(\sum_{i=0}^{T-1}\gamma_i\Delta
W\left(i\right)\right)^2=\frac{1}{T^2}\sum_{i=0}^{T-1} \Ex_\vartheta \gamma
_i^2=\frac{\Ex_\vartheta \gamma_0^2}{T}\leq \frac{K}{T},\\
J_3&=\frac{1}{T^2}\Ex_\vartheta \left(\sum_{i=0}^{T-1}\left[\Delta W\left(i\right)^2-1
  \right]  \right)^2=\frac{1}{T^2}\sum_{i=0}^{T-1} \Ex_\vartheta \left[\Delta W\left(i\right)^2-1
  \right]^2\leq \frac{K}{T}.
\end{align*}
The second proposition of the  Lemma 2 follows from these estimates and
\eqref{2.6}.
\begin{lemma}
\label{L3}
The function $\Phi \left(x\right)$ (see \eqref{fi}) has the following properties
\begin{align}
\label{2.8}
&\lim_{x\rightarrow  0+}\Phi \left(x\right)=\frac{1}{2},\\
&\lim_{x\rightarrow  \infty }\Phi \left(x\right)=0,
\label{2.9}\\
&\Phi '\left(x\right)<0, \quad {\rm for}\quad x>0.
\label{2.10}
\end{align}

\end{lemma}
{\bf Proof.} From \eqref{2.2} we obtain the representations
\begin{align*}
\Phi
\left(x\right)&=\frac{1}{2}-\frac{x}{3!}+\frac{x^2}{4!}-\frac{x^3}{5!}+\ldots\\ 
\Phi'\left(x\right)&=-\left(\frac{1}{3!}-\frac{2x}{4!}\right)-
\left(\frac{3x^2}{5!}-\frac{4x^3}{6!}\right)-\ldots
\end{align*}
which allow to verify  the limits \eqref{2.8} and  \eqref{2.9} and as well the
estimate \eqref{2.10} for $x<2$. For $x\geq 2$ this estimate follows from
the explicit expression for this derivative
\begin{align*}
\Phi
'\left(x\right)=\frac{1}{x^2}\left(\frac{2}{x}-1\right)-\left(\frac{2}{x^3}+\frac{1}{x^2}\right)
e^{-x} .
\end{align*}

Let us consider the equation \eqref{2.11} for $\alpha _T$, here $\zeta _T$  and
$\eta _T$ are defined in \eqref{2.2} and \eqref{2.12a} respectively.
 Due to Lemma \ref{L3} this equation has not more than one
solution. Recall that  ${\cal A}_T$ is the following  event: the equation \eqref{2.11} has
solution and consider the statistics $\beta _T$  defined in \eqref{2.12} (here
$c_0,c_1$ are the constants from the condition \eqref{1.6}).

\begin{lemma}
\label{L4}
 Under the condition \eqref{1.6} the estimate $\beta _T$ is
$\sqrt{T}$-consistent for $\lambda +\mu $. Moreover, for some constant $C>0$
which does not depend on $T$ and $\vartheta $ we have the property
\begin{equation}
\label{2.13}
\Ex_\vartheta \left[\sqrt{T}\left(\beta _T-\left(\lambda +\mu \right)\right)\right]^2<C.
\end{equation}
\end{lemma}
{\bf Proof.} It follows from Lemmae \ref{L1} and \ref{L2} that 
\begin{equation}
\label{2.14}
\zeta _T=\bar Y^2+2D\Phi \left(\lambda +\mu \right)+ \varepsilon _1\left(T\right).
\end{equation}
Here and below we have for $\varepsilon _i\left(T\right), i=1,2\ldots $  the
 estimates 
\begin{align*}
\Ex_\vartheta\left(\sqrt{T}\varepsilon _i\left(T\right)\right)^2<C.
\end{align*}

By Lemma \ref{L1}, estimates \eqref{1.8}, \eqref{1.9} and the boundedness   of
$\Phi \left(x\right)$ we obtain as well the relation
\begin{equation}
\label{2.15}
\zeta _T=\bar Y^2+2\tilde \eta  _T\Phi \left(\lambda +\mu \right)+ \varepsilon_2\left(T\right).
\end{equation}

If we have the event ${\cal A}_T$ then  it follows from \eqref{2.11} and \eqref{2.15} that
\begin{align*}
2\tilde\eta  _T\Phi \left(\alpha _T\right)=2\tilde\eta _T \Phi
\left(\lambda +\mu \right)+\varepsilon _3\left(T\right). 
\end{align*}
This relation, Lemma \ref{L3} and the separation from zero by a positive
constant of the estimator $\tilde\eta  _T$ (see Corollary to Lemma
\ref{L1}) yield for $\omega \in {\cal A}_T$
\begin{align*}
\Phi \left(\alpha _T\right)-\Phi \left(\lambda +\mu \right)=\varepsilon
_4\left(T\right). 
\end{align*}
Therefore from Lemma \ref{L3} we obtain
\begin{equation}
\label{2.16}
\Ex_\vartheta \left\{\1_{\left\{A_T\right\}} \sqrt{T}\left(\beta _T-\left(\lambda +\mu
\right) \right)\right\}^2<C. 
\end{equation}
If $\omega \in {\cal A}_T^c$ then the equation
\begin{equation}
\label{2.17}
\bar Y^2+2D\Phi \left(\lambda +\mu \right)+\gamma
_3\left(T\right)=\left(\frac{X\left(T\right)}{T}\right)^2+2\tilde\eta  _T\Phi \left(x\right) 
\end{equation}
has no solution $x\in \left[2c_0,2c_1\right]$.

It follows from \eqref{2.17}, Lemma \ref{L1} and the Corollary that the
equation
\begin{align*}
\Phi \left(x\right)=\Phi \left(\lambda +\mu \right)+\varepsilon  _4\left(T\right)
\end{align*}
has no solution for $x\in \left[2c_0,2c_1\right]$.

 Hence we can write ${\cal A}_T^c\subset \left\{\left|\varepsilon
 _4\left(T\right)\right|>\alpha \right\}={\cal B}_T$ for some positive constant $\alpha $
 which does not depend on $T$. This allow us to write
\begin{align*}
\Pb\left({\cal A}_T^c\right)\leq \Pb\left({\cal B}_T\right)=\Pb\left\{\left|\varepsilon
_4\left(T\right)\right|>\alpha \right\} &\leq \Pb\left\{\left|\sqrt{T}\varepsilon
_4\left(T\right)\right|^2>\alpha ^2T\right\}\\
& \leq \frac{\Ex_\vartheta \left|\sqrt{T}\varepsilon
_4\left(T\right)\right|^2}{\alpha ^2T}<\frac{C}{T}.
\end{align*} 
This estimate and \eqref{2.16} prove the Lemma \ref{L4}.

\bigskip
 
{\bf Proof of the Theorem}  \ref{T1}. 
The obtained results allow us to prove that the estimators defined in
\eqref{2.19} are $\sqrt{T}$-consistent. 
Indeed, from the obvious equality
\begin{align*}
\hat\lambda _T-\lambda =\beta _T\frac{\frac{X\left(T\right)}{T}-y_1}{y_2-y_1}
-\beta _T\frac{\lambda }{\lambda +\mu } +\frac{\lambda }{\lambda +\mu
}\left(\beta _T-\left(\lambda +\mu \right)\right) 
\end{align*}
and \eqref{1.7} we obtain by Lemma \ref{L4} the estimate
\begin{align*}
\Ex_\vartheta \left(\sqrt{T}\left(\hat\lambda _T-\lambda \right)\right)^2&\leq 2\Ex_\vartheta\left[
   \sqrt{T}  \beta
   _T\left(\frac{\frac{X\left(T\right)}{T}-y_1}{y_2-y_1}-\frac{\lambda
   }{\lambda +\mu } \right) \right]^2\\
&\quad +2\left(\frac{\lambda }{\lambda +\mu }\right)^2\Ex_\vartheta \left[
  \sqrt{T}\left(\beta _T-\left(\lambda +\mu  \right)\right)\right]^2<C.
\end{align*}
The second inequality in \eqref{2.19} can be proved by the same way. 

Therefore the estimator $\hat\vartheta _T=\left(\hat\lambda _T,\hat\mu
_T\right)$ is $\sqrt{T}$-consistent.

\section{One-step MLE}

Our goal is to construct the asymptotically efficient estimator-process of the
parameter $\vartheta =\left(\lambda ,\mu \right)\in \Theta $. We do it in two
steps.  First we obtain by the observations $X^{T^\delta }=\left(X_t,0\leq
t\leq T^\delta \right)$ on the learning interval $\left[0,T^\delta \right]$
the  method of moments estimator $ \hat\vartheta _{T^\delta
}=(\hat\lambda _{T^\delta } ,\hat\mu _{T^\delta } )$ studied in the preceding
section.  Here $\delta \in \left(\frac{1}{2},1\right)$. This estimator by
Theorem \ref{T1} satisfies the condition:
\begin{align*}
\sup_{\vartheta \in K}T^{{\delta}{} }\Ex_\vartheta \left|\hat\vartheta
_{T^\delta }-\vartheta \right|^2\leq C,
\end{align*}
here the constant $C>0$ does not depend on $T$ and $\vartheta \in \Theta $.
Remind that $\Theta
=\left(c_0,c_1\right)\times\left(c_0,c_1\right)$. Introduce the additional
condition
\bigskip

${\cal M}\left(N\right).$ {\it For some} $N\geq 2$
\begin{equation}
\label{34}
\frac{c_0}{\left(y_1-y_2\right)^2}>\frac{2N+9}{4}.
\end{equation}

Having this preliminary estimator $\hat\vartheta
_{T^\delta } $ we propose one-step MLE which is based on
one modification of the score-function
\begin{align*}
\Delta_t \left(\vartheta ,X^t\right)=\frac{1}{\sqrt{t}}\int_{0}^{t}{\dot
  m(\vartheta ,s)}\,\left[ {\rm d}X_s- m(\vartheta
  ,s){\rm d}s\right],\qquad T^\delta \leq t\leq T
\end{align*}
 as follows
\begin{align}
\label{os}
\vartheta _{t,T}^\star=\hat\vartheta _{T^\delta }+t^{-1} \II_t(\hat\vartheta
_{T^\delta } )^{-1}  \int_{T^\delta }^{t}{\dot 
  m(\hat\vartheta _{T^\delta } ,s)}\,\left[ {\rm d}X_s- m(\hat\vartheta _{T^\delta }
  ,s){\rm d}s\right].
\end{align}
Here   the vector 
\begin{align*}
\dot   m(\vartheta ,s)=\left(y_1-y_2\right)\frac{\partial \pi_\lambda
    \left(s,\vartheta \right)}{\partial  \vartheta }=\left(y_1-y_2\right)\left(\frac{\partial \pi
  \left(t,\vartheta \right)}{\partial \lambda },\frac{\partial \pi
  \left(t,\vartheta \right)}{\partial \mu  } \right)^* 
\end{align*}
and the empirical   Fisher information matrix $\II_t(\vartheta ) $ is
\begin{align*}
\II_t(\vartheta )=\frac{1}{t}\int_{T^\delta }^{t}{\dot   m(\vartheta ,s)}{\dot
  m(\vartheta ,s)}^*{\rm d}s\longrightarrow  \II(\vartheta )
\end{align*}
as $t\rightarrow \infty $ by the law of large numbers. 
Here $\II(\vartheta ) $  is the Fisher information matrix
\begin{align*}
\II(\vartheta )=\left(y_1-y_2\right)^2\Ex_\vartheta \frac{\partial \pi
    \left(s,\vartheta \right)}{\partial  \vartheta }\frac{\partial \pi
    \left(s,\vartheta \right)^*}{\partial  \vartheta }.
\end{align*}

Let us change the variable $\tau =tT^{-1}\in \left[0,1\right]$ and introduce the random process
 $\vartheta
_{ T}^\star\left(\tau \right), \tau _\delta\leq \tau \leq 1 ,$ here $\vartheta
_{ T}^\star\left(\tau \right)=\vartheta
_{ \tau T,T}^\star$ and   $\tau _\delta =T^{\delta -1}\rightarrow 0$.

\begin{theorem}
\label{T2} Suppose that $\vartheta \in \Theta $, $\delta \in (\frac{1}{2},
1)$ and the condition ${\cal M}\left(2\right)$ holds, then the one-step
MLE-process is consistent: for any $\nu >0$ and any $\tau \in (0,1]$
\begin{equation}
\label{35a}
\Pb_{\vartheta _0}\left\{ \left| \vartheta
_T^\star\left(\tau \right)-\vartheta _0 \right|>\nu \right\}\rightarrow 0
\end{equation}
and it is asymptotically normal
\begin{equation}
\label{36}
\sqrt{\tau T}\left(\vartheta _T^\star\left(\tau \right)-\vartheta _0\right)\Longrightarrow {\cal
  N}\left(0,\II(\vartheta_0 )^{-1} \right).
\end{equation}

\end{theorem}
{\bf Proof.} Let us denote 
   the   partial derivatives
\begin{align*}
\dot\pi_\lambda  \left(t,\vartheta \right)=\frac{\partial \pi
  \left(t,\vartheta \right)}{\partial \lambda },\qquad \dot\pi_\mu 
\left(t,\vartheta \right)=\frac{\partial \pi
  \left(t,\vartheta \right)}{\partial \mu },\quad \ddot\pi_{\lambda ,\lambda}
\left(t,\vartheta \right)=\frac{\partial^2 \pi 
  \left(t,\vartheta \right)}{\partial \lambda^2 },
\end{align*}
and so on.  
\begin{lemma}
\label{L5} Suppose that $\vartheta \in \Theta $ and $N>1$. 
If the condition 
\begin{equation}
\label{c1}
\frac{c_0}{ \left(y_1-y_2\right)^2}>\frac{N+1}{4}
\end{equation}
 holds, then 
\begin{align}
\label{36b}
\sup_{\vartheta \in \Theta } \Ex_{\vartheta _0}\left(\left|\dot\pi_\lambda
\left(t,\vartheta \right) \right|^N+ \left|\dot\pi_\mu 
\left(t,\vartheta \right) \right|^N \right) <C_1,
\end{align}
and if the condition 
\begin{equation}
\label{c2}
\frac{c_0}{ \left(y_1-y_2\right)^2}>\frac{2N+9}{4}
\end{equation}
 holds, then 
\begin{align}
\label{36c}
\sup_{\vartheta \in \Theta } \Ex_{\vartheta _0}\left(\left|\ddot\pi_{\lambda,\lambda }
\left(t,\vartheta \right) \right|^N+ \left|\ddot\pi_{\lambda,\mu  }
\left(t,\vartheta \right) \right|^N+ \left|\ddot\pi_{\mu ,\mu  }
\left(t,\vartheta \right) \right|^N \right) <C_2.
\end{align}
Here the constants $C_1>0,C_2>0$ do not depend on $t$.
\end{lemma}
{\bf Proof.} For simplicity of exposition we write
\begin{align*}
\dot \pi _\lambda \left(t,\vartheta \right)=\dot
\pi _\lambda ,\qquad  \dot \pi _\mu  \left(t,\vartheta \right)=\dot
\pi _\mu ,\qquad   \pi \left(t,\vartheta \right)=\pi.
\end{align*}
  By the formal differentiation  of
\begin{align}
{\rm d}\pi &=\left[\mu -\left(\lambda +\mu \right)\pi
   -\pi \left(1-\pi
\right)
\left(y_1-y_2\right)\left(y_2+\left(y_1-y_2\right)\pi \right)\right]{\rm d}t\nonumber\\
&\quad +\pi\left(1-\pi
\right)\left(y_1-y_2\right)\;{\rm d}X_t.
\label{37}
\end{align}
 we obtain the equations
\begin{align}
\label{04}
{\rm d}\dot \pi _\lambda =&-\pi \,
{\rm d}t - \dot \pi _\lambda \left[ \lambda +\mu + \left(1-2\pi
  \right)\left(y_1-y_2\right)\left[y_2+\left(y_1-y_2\right)\pi \right]\right.\nonumber\\
&\left. +\pi \left(1-\pi
  \right)\left(y_1-y_2\right)^2 \right]
{\rm d}t +\dot \pi _\lambda \left(1-2\pi \right)\left(y_1-y_2\right)\,{\rm
  d}X\left(t\right),\\
{\rm d}\dot \pi _\mu  =&\left[1-\pi\right] 
{\rm d}t -\dot \pi _\mu \left[ \lambda +\mu + \left(1-2\pi
  \right)\left(y_1-y_2\right)\left[y_2+\left(y_1-y_2\right)\pi \right]\right.\nonumber\\
&\left. +\pi \left(1-\pi
  \right)\left(y_1-y_2\right)^2 \right]
{\rm d}t +\dot \pi _\mu \left(1-2\pi \right)\left(y_1-y_2\right)\,{\rm
  d}X\left(t\right).
\label{05}
\end{align}
If we denote the true value of the parameters by $ \vartheta _0$ and $\pi
\left(t,\vartheta_0 \right)=\pi^o$ etc., then these equations for $\vartheta
=\vartheta _0$ become
\begin{align}
\label{06}
{\rm d}\dot \pi _\lambda^o =&-\pi^o \,
{\rm d}t - \dot \pi _\lambda^o \left[ \lambda_0 +\mu_0 +\pi^o \left(1-\pi^o
  \right)\left(y_1-y_2\right)^2 \right]
{\rm d}t\nonumber\\
& +\dot \pi _\lambda^o \left(1-2\pi^o \right)\left(y_1-y_2\right)\,{\rm
  d}\bar W\left(t\right),\\
{\rm d}\dot \pi _\mu^o  =&\left[1-\pi^o\right] 
{\rm d}t - \dot \pi _\mu^o \left[ \lambda_0 +\mu_0 +\pi  ^o\left(1-\pi^o
  \right)\left(y_1-y_2\right)^2 \right]
{\rm d}t\nonumber\\
& +\dot \pi _\mu ^o\left(1-2\pi \right)\left(y_1-y_2\right)\,{\rm
  d}\bar W\left(t\right),
\label{07}
\end{align}
here
\begin{align}
\label{08}
{\rm d}\pi ^o=\left[\mu _0-\left(\lambda _0+\mu _0\right)\pi ^o\right]{\rm
  d}t+\pi ^o\left(1-\pi ^o\right) \left(y_1-y_2\right)\,{\rm
  d}\bar W\left(t\right).
\end{align}

This system of  linear for $\dot \pi _\lambda $ and $\dot \pi _\mu  $   equations can be re-written as follows ($x_t=\pi
^o,y_t=\dot \pi_\lambda  ^o, z_t=\dot \pi_\mu   ^o ,a=\lambda_0 +\mu_0 , b=y_1-y_2$)
\begin{align}
\label{xt}
{\rm d}x_t &=\left[\mu _0-ax_t\right]{\rm d}t+bx_t\left(1-x_t\right){\rm
  d}\bar W_t,\\
\label{yt}
{\rm d}y_t&=-x_t{\rm d}t-\left[a+b^2x_t\left(1-x_t\right)\right]y_t{\rm
  d}t+b\left(1-2x_t\right)y_t\;{\rm d}\bar W_t,\\
{\rm d}z_t&=\left[1-x_t\right]{\rm d}t-\left[a+b^2x_t\left(1-x_t\right)\right]z_t{\rm
  d}t+b\left(1-2x_t\right)z_t\;{\rm d}\bar W_t.
\label{zt}
\end{align}

Note that as $\lambda _0>0$ and $\mu _0>0$ the process $\pi\left(t,\vartheta
_0\right)=x_t\in \left(0,1\right) $ is ergodic with two reflecting borders 0 and
1. Therefore the process $\pi ^o\left(t,\vartheta _0\right)$ is ergodic with
 the invariant density 
\begin{align*}
f\left(\vartheta _0,x\right)&= 
\frac{\left[x\left(1-x\right)\right]^{\frac{2\left(\mu_0 -\lambda_0 \right)}{
      \left(y_1-y_2\right)^2}-2}}{G\left(\vartheta_0 \right)}\;
\exp\left\{-\frac{2\mu_0 +2\left(\lambda_0 -\mu_0 \right)x}{\left(y_1-y_2\right)^2
  x\left(1-x\right)}\right\} \nonumber\\
&=\frac{\left[x\left(1-x\right)\right]^{\gamma \left(\mu_0 -\lambda_0 \right)-2}}{G\left(\vartheta_0 \right)}\;
\exp\left\{-\frac{\gamma \mu_0 }{x}-\frac{\gamma \lambda_0 }{1-x} \right\}
\end{align*}
here we denoted $\gamma =2\left(y_1-y_2\right)^{-2}$ and
$G\left(\vartheta _0\right) $ is the normalizing constant
\begin{align*}
G\left(\vartheta _0\right) =\int_{0}^{1}{\left[x\left(1-x\right)\right]^{\gamma \left(\mu_0 -\lambda_0 \right)-2}}{}\;
\exp\left\{-\frac{\gamma \mu_0 }{x}-\frac{\gamma \lambda_0 }{1-x} \right\} \;{\rm d}x.
\end{align*}

The processes $y_t$ and $z_t$ have explicite expressions
\begin{align}
y_t&=-\int_{0}^{t}\exp\left\{-\int_{v}^{t}
\left[a+b^2x_s\left(1-x_s\right)-\frac{b^2}{2}
  \left(1-2x_s\right)^2\right]{\rm d}s\right.\nonumber\\
\label{09}
&\qquad \qquad \qquad \left.+b\int_{v}^{t}
\left(1-2x_s\right){\rm d}\bar W_s\right\}x_v\,{\rm d}v,\\
z_t&=\int_{0}^{t}\exp\left\{-\int_{v}^{t}
\left[a+b^2x_s\left(1-x_s\right)-\frac{b^2}{2}
  \left(1-2x_s\right)^2\right]{\rm d}s\right.\nonumber\\
&\qquad \qquad \qquad \left.+b\int_{v}^{t}
\left(1-2x_s\right){\rm d}\bar W_s\right\}\left[1-  x_v\right]\,{\rm d}v.
\label{10}
\end{align}
Let us put ${\rm x}_s=\frac{1}{2} -{x}_s$. Then  we have 
\begin{align*}
x_s\left(1-x_s\right)-\frac{1}{2}  \left(1-2x_s\right)^2=-3{\rm x}_s^2+\frac{1}{4}
\end{align*}
and
\begin{align*}
y_t&=\int_{0}^{t} \bigl( {\rm x}_v-\frac{1}{2}
\bigr)e^{-\left(a+\frac{b^2}{4}\right)\left(t-v\right)}
\exp\left\{{3b^2\int_{v}^{t}{\rm x}_s^2{\rm d}s +2b\int_{v}^{t}{\rm x}_s{\rm
    d}\bar W_s}\right\}\;{\rm d}v.
  \end{align*}

To estimate the moments $\Ex_{\vartheta_0}  \left|y_t\right|^N$ we note that
$\left|{\rm x}_v-\frac{1}{2}\right| \leq \frac{1}{2}$ and use the H\"older
inequality 
\begin{align*}
\left(\int_{0}^{t}\left|f\left(v\right)g\left(v\right)\right|{\rm
  d}v\right)^N\leq \left(\int_{0}^{t}\left|f\left(v\right)\right|^{\frac{N}{N-1}}{\rm
  d}v\right)^{N-1} \int_{0}^{t}\left|g\left(v\right)\right|^N{\rm d}v
\end{align*}
with $f\left(v\right)=\exp\left\{-a\left(t-v\right)\varepsilon \right\}$ and 
\begin{align*}
g\left(v\right)=\exp\left\{-\left(a\left(1-\varepsilon
\right)+\frac{b^2}{4}\right)\left(t-v\right) {+3b^2\int_{v}^{t}{\rm x}_s^2{\rm
    d}s +2b\int_{v}^{t}{\rm x}_s{\rm     d}\bar W_s}   \right\}, 
\end{align*}
here $\varepsilon >0$. This yields the estimate
\begin{align*}
\Ex_{\vartheta_0}  \left|y_t\right|^N\leq C\left(N,\varepsilon \right)
\int_{0}^{t} e^{-N\left(a\left(1-\varepsilon
  \right)+\frac{b^2}{4}\right)\left(t-v\right)}\Ex_{\vartheta_0}e^{
  {{3Nb^2}\int_{v}^{t}{\rm x}_s^2{\rm d}s +2Nb\int_{v}^{t}{\rm x}_s{\rm 
    d}\bar W_s}   }{\rm d}v,
\end{align*}
here the constant $C\left(N,\varepsilon \right) >0 $ does not depend on $t$. 
Further, we can write
\begin{align*}
&\Ex_{\vartheta_0}\exp\left\{
  3Nb^2\int_{v}^{t}{\rm x}_s^2{\rm d}s +2Nb\int_{v}^{t}{\rm x}_s{\rm 
    d}\bar W_s\right\}\\
&\qquad  =\Ex_{\vartheta_0}\left(\exp\left\{{2Nb}\int_{v}^{t}{\rm x}_s{\rm 
    d}\bar W_s- {2N^2b^2}\int_{v}^{t}{\rm x}_s^2{\rm d}s\right\}\right.\\
&\qquad \qquad \qquad \left.\exp\left\{Nb^2\left(2N+3\right)\int_{v}^{t}{\rm x}_s^2{\rm d}s\right\}\right)\\
&\qquad  \leq \exp\left\{ \frac{N b^2}{4}\left({2N}+3\right) \left(t-v\right)\right\}
\end{align*}
because ${\rm x}_s^2\leq 1/4  $ and 
\begin{align*}
\Ex_{\vartheta_0}\exp\left\{{{2Nb}\int_{v}^{t}{\rm x}_s{\rm 
    d}\bar W_s- {2N^2b^2}\int_{v}^{t}{\rm x}_s^2{\rm d}s}\right\}=1.
\end{align*}
Therefore, 
\begin{align*}
\Ex_{\vartheta_0}  \left|y_t\right|^N&\leq C\left(N,\varepsilon \right)
\int_{0}^{t} e^{-{N}\left(a\left(1-\varepsilon
  \right)+\frac{b^2}{4}-\frac{b^2}{4}\left(
  2N+3\right)\right)\left(t-v\right)}{\rm d}v\\
&=C\left(N,\varepsilon \right)
\int_{0}^{t} e^{-N\left(a\left(1-\varepsilon \right)
  -\frac{b^2}{2}\left(N+1\right)\right)\left(t-v\right)}{\rm
  d}v.  
\end{align*}
We see that if 
\begin{align*}
\frac{\lambda _0+\mu _0}{\left(y_1-y_2\right)^2 }>\frac{1}{2}+\frac{N}{2},
\end{align*}
then $\Ex_{\vartheta _0}  \left|y_t\right|^N\leq C  $. In particular,  
if in the condition \eqref{c1} we put $N=2$ and choose sufficiently
small $\varepsilon >0$, then  we obtain the estimate
\begin{equation}
\label{46}
\sup_{\vartheta_0 \in \Theta }\Ex_{\vartheta_0}  \left|\frac{\partial \pi
  \left(t,\vartheta_0 \right)}{\partial \lambda }\right|^2\leq C,
\end{equation}
here the constant $C>0$ does not depend on $t$.

We need as well to estimate  the derivatives \eqref{04}, \eqref{05}   for
the values $\vartheta \not=\vartheta _0$. 
The equation for $\dot \pi _\lambda $ becomes 
\begin{align}
{\rm d}\dot \pi _\lambda =&-\pi \,
{\rm d}t - \dot \pi _\lambda \left[ \lambda +\mu + \left(1-2\pi
  \right)\left(y_1-y_2\right)^2\left(\pi-\pi ^0\right)\right.\nonumber\\
&\left. +\pi \left(1-\pi
  \right)\left(y_1-y_2\right)^2 \right]
{\rm d}t +\dot \pi _\lambda \left(1-2\pi \right)\left(y_1-y_2\right)\,{\rm
  d}\bar W\left(t\right).
\label{47}
\end{align}
Hence if we put $a=\lambda +\mu, y_t=\dot \pi _\lambda $ and
$b=y_1-y_2$, then we obtain the equation
\begin{align*}
{\rm d}y_t&=-x_t{\rm
  d}t-\left[a+b^2\left(1-2x_t\right)\left(x_t-x_t^0\right)+b^2x_t\left(1-x_t\right)\right]y_t{\rm
  d}t \\ 
& \qquad +b\left(1-2x_t\right)y_t\,{\rm d}\bar W_t.
\end{align*}
The solution of this equation can be written explicitly like \eqref{09} but
with additional term  $b^2\left(1-2x_t\right)\left(x_t-x_t^0\right) $ in the
exponent. This term   satisfies the inequality
\begin{align*}
\left(1-2x_t\right)\left(x_t-x_t^0\right)\geq -1.
\end{align*}
Hence if we repeat the evaluation of the $\Ex_{\vartheta
  _0}\left|y_t\right|^2$ as it was done above, then for it boundness we obtain the condition
\begin{align*}
\frac{\lambda +\mu }{\left(y_1-y_2\right)^2}>\frac{3}{2}+\frac{N}{2}.
\end{align*}

For the second derivative $\ddot \pi= \ddot \pi_{\lambda ,\lambda }\left(t,\vartheta \right)$
 we obtain  the similar estimates of the moments as follows. The equation  for $\ddot \pi $ is
\begin{align*}
{\rm d}\ddot\pi&=-y_t\left[
  2-2b^2y_t\left(x_t-x_t^0\right)+2b^2y_t\left(1-2x_t\right)\right]{\rm d}t-2by_t^2\,{\rm
  d}\bar W_t\\
&\quad -\ddot \pi
\left[a+b^2\left(1-2x_t\right)\left(x_t-x_t^0\right)+b^2x_t\left(1-x_t\right)
  \right] {\rm d}t+b\ddot \pi \left(1-2x_t\right){\rm d}\bar W_t.
\end{align*}
Let us write it as 
\begin{align*}
{\rm d}\ddot\pi&=A\left(t\right){\rm d}t+B\left(t\right)\,{\rm
  d}\bar W_t-\ddot\pi \left[a+C\left(t\right)\right]{\rm d}t+\ddot\pi_t D\left(t\right)\,{\rm
  d}\bar W_t
\end{align*}
in obvious notations. Hence the solution of it is
\begin{align*}
\frac{\partial^2 \pi
  \left(t,\vartheta \right)}{\partial \lambda^2 }=\int_{0}^{t}
e^{-\int_{v}^{t}\left[a+C\left(s\right)-\frac{1}{2}D\left(s\right)^2\right]{\rm
    d}s +\int_{v}^{t}D\left(s\right){\rm d}\bar W_s }\left[A\left(v\right){\rm
    d}v+B\left(v\right){\rm d}\bar W_v \right] .
\end{align*}
We have the corresponding estimate 
\begin{align*}
&C\left(s\right)-\frac{D\left(s\right)^2}{2}=
b^2\left(1-2x_s\right)\left(x_s-x_s^0\right)+\frac{b^2}{2} \left[2x_s\left(1-x_s\right)-\left(1-2x_s\right) 
^2 \right]\\
&\qquad \quad \geq
-b^2-3b^2\left(x-\frac{1}{2}\right)^2+\frac{b^2}{4}=-\frac{3b^2}{4}-3b^2\left(x-\frac{1}{2}\right)^2\geq  -\frac{3b^2}{2}
\end{align*}
because $\left(x-\frac{1}{2}\right)^2\leq \frac{1}{4}$. Therefore,
\begin{align*}
a-\frac{3}{2}b^2-\frac{2N+3}{4}b^2=a-\frac{9}{4}b^2-\frac{N}{2}b^2
\end{align*}
and  if
\begin{align*}
\frac{\lambda +\mu }{\left(y_1-y_2\right) ^2}>\frac{9}{4}+\frac{N}{2},
\end{align*}
then we obtain
\begin{align*}
\Ex_{\vartheta _0}\left|\frac{\partial^2 \pi
  \left(t,\vartheta \right)}{\partial \lambda^2 }  \right|^N<C.
\end{align*}
Hence under condition \eqref{c2} we have
\begin{align*}
\sup_{\vartheta \in \Theta }\Ex_{\vartheta _0}\left|\frac{\partial^2 \pi
  \left(t,\vartheta \right)}{\partial \lambda^2 }  \right|^N<C.
\end{align*}
The similar estimates can be obtained for the other derivatives. Lemma
\ref{L5} is proven.
\begin{lemma}
\label{L6}
The solutions $\left(x_t,y_t,z_t\right)$
 of the equations \eqref{xt}-\eqref{zt} have ergodic properties. In
particular, we have the following  mean square convergence 
\begin{align*}
\frac{1}{T}\int_{0}^{T}\dot m_\lambda \left(t,\vartheta_0 \right)^2{\rm
  d}t&=\frac{ b^2}{T}\int_{0}^{T}y_t^2\;{\rm
  d}t\longrightarrow {\rm I}_{11}\left(\vartheta_0 \right), \\
\frac{1}{T}\int_{0}^{T}\dot m_\lambda \left(t,\vartheta_0 \right)\dot m_\mu
  \left(t,\vartheta_0 \right){\rm 
  d}t&=\frac{ b^2}{T}\int_{0}^{T}y_tz_t\;{\rm
  d}t\longrightarrow {\rm I}_{12}\left(\vartheta_0 \right),\\
\frac{1}{T}\int_{0}^{T}\dot m_\mu  \left(t,\vartheta_0 \right)^2\;{\rm
  d}t&=\frac{ b^2}{T}\int_{0}^{T}z_t^2\;{\rm
  d}t\longrightarrow {\rm I}_{22}\left(\vartheta_0 \right),
\end{align*}
\end{lemma}
{\bf Proof.} For the proof of the invariant measure  existence see
\cite{PCh09}, section 4.2. Note that the equations \eqref{xt}-\eqref{zt} do
not coincide with that of \cite{PCh09}, because there it is supposed that
$\lambda =\mu $ and $y_1=1,y_2=0$, but the arguments given there are
directly applied to the system of equations \eqref{xt}-\eqref{zt} too.

Recall that the strong mixing coefficient $\alpha \left(t\right)$ for ergodic diffusion process
\eqref{xt}  satisfies the estimate
\begin{align*}
\alpha \left(t \right)<e^{-c\left|t\right|}.
\end{align*}
For the proof see Theorem in \cite{V87}. To check the conditions of this
theorem we change the variables in the equation \eqref{xt}
\begin{align*}
\xi _t= g\left(x_t\right),\qquad g\left(x\right)= \int_{1/2}^{x}\frac{{\rm
    d}v}{bv\left(1-v\right)},\quad x\in \left(0,1\right)
\end{align*}
and obtain the stochastic differential equation
\begin{align*}
{\rm d}\xi _t=A\left(\xi _t\right){\rm d}t+{\rm d}W_t,\qquad \xi _0= g\left(x_0\right),\qquad 0\leq t\leq T.
\end{align*}
The process $\xi _t,t\geq 0$ has ergodic properties  and the drift coefficient
$A\left(\cdot \right)$ satisfies  the conditions of this theorem. 

 Now to verify the convergence
\begin{align}
&\Ex_{\vartheta _0}\left(\frac{ 1}{T}\int_{0}^{T}y_t^2\;{\rm  d}t-\frac{ 1}{T}\int_{0}^{T}\Ex_{\vartheta _0}y^2
_t\;{\rm  d}t\right)^2 \nonumber\\
&\qquad \qquad =\Ex_{\vartheta _0}\left(\frac{ 1}{T}\int_{0}^{T}\left[y_t^2-\Ex_{\vartheta _0}y^2
_t\right]\;{\rm  d}t\right)^2\longrightarrow  0
\label{qq}
\end{align}
we can apply the result of the following lemma. 

\begin{lemma}
\label{L7}
Let $\left\{Y_t,t>0\right\}$ be a stochastic process with zero  mean and for
some $m>2$ and $k\geq 1$
\begin{align*}
\Ex \left|Y_t \right|^{m\left(2k-1\right)}<C_1,\qquad \int_{0}^{\infty
}t^{k-1}\left[ \alpha \left(t\right)\right]^{\left(m-2\right)/m}\;{\rm d}t<C_2,
\end{align*}
here $\alpha \left(t\right)$ is  the strong mixing the coefficient.
Then
\begin{equation*}
\Ex\left|\int_{0}^{T}Y_t\;{\rm d}t \right|^{2k}\leq C_3\;T^k.
\end{equation*}
\end{lemma}
{\bf Proof}. For proof see   Lemma 2.1 in \cite{K66}.

Therefore if we put $Y_y=y_t^2-\Ex_{\vartheta _0}y^2_t$, $m=3$ and $k=1$, then we obtain the
convergence \eqref{qq}.

Let us verify the consistency of the one-step MLE-process.
We can write
\begin{align*}
&\Pb_{\vartheta _0}\left\{ \left| \vartheta
_T^\star\left(\tau \right)-\vartheta _0 \right|>\nu \right\}\leq  \Pb_{\vartheta _0}\left\{ \left| \hat\vartheta
_{T^\delta }-\vartheta _0 \right|>\frac{\nu }{2}\right\}\\
&\quad +\Pb_{\vartheta _0}\left\{\left| \frac{\II_{\tau T}(\hat\vartheta
_{T^\delta } )^{-1}}{\tau T}  \int_{T^\delta }^{\tau T}{\dot 
  m(\hat\vartheta _{T^\delta } ,s)}\,\left[ {\rm d}X_s- m(\hat\vartheta _{T^\delta }
  ,s){\rm d}s\right]\right|>\frac{\nu }{2}\right\}.
\end{align*}
For the first probability by the Theorem \ref{T1} we have
\begin{align*}
\Pb_{\vartheta _0}\left\{ \left| \hat\vartheta
_{T^\delta }-\vartheta _0 \right|>\frac{\nu }{2}\right\}\leq \frac{4}{\nu ^2}\Ex_{\vartheta _0}\left| \hat\vartheta
_{T^\delta }-\vartheta _0 \right|^2 \leq \frac{C}{\nu ^2T^\delta }\rightarrow 0.
\end{align*}
 The second probability  can be  evaluated as follows
\begin{align*}
&\Pb_{\vartheta _0}\left\{\left| \frac{\II_{\tau T}(\hat\vartheta
_{T^\delta } )^{-1}}{\tau T}  \int_{T^\delta }^{\tau T}{\dot 
  m(\hat\vartheta _{T^\delta } ,s)}\,\left[ {\rm d}X_s- m(\hat\vartheta _{T^\delta }
  ,s){\rm d}s\right]\right|>\frac{\nu }{2}\right\}\\
&\; \leq  \Pb_{\vartheta _0}\left\{\left| \frac{\II_{\tau T}(\hat\vartheta
_{T^\delta } )^{-1}}{\tau T}  \int_{T^\delta }^{\tau T}{\dot 
  m(\hat\vartheta _{T^\delta } ,s)}\, {\rm d}\bar W_s\right|>\frac{\nu }{4}\right\}\\
&\qquad  +\Pb_{\vartheta _0}\left\{\left| \frac{\II_{\tau T}(\hat\vartheta
_{T^\delta } )^{-1}}{\tau T}  \int_{T^\delta }^{\tau T}{\dot 
  m(\hat\vartheta _{T^\delta } ,s)}\,\Delta m\left(\hat\vartheta _{T^\delta },s\right){\rm d}s\right|>\frac{\nu }{4}\right\},
\end{align*}
here $ \Delta m\left(\hat\vartheta _{T^\delta },s\right)=m(\vartheta _0
  ,s)- m(\hat\vartheta _{T^\delta }
  ,s)$.  
We can write
\begin{align*}
&\left| \frac{\II_{\tau T}(\hat\vartheta
_{T^\delta } )^{-1}}{\tau T}  \int_{T^\delta }^{\tau T}{\dot 
  m(\hat\vartheta _{T^\delta } ,s)}\, {\rm d}\bar W_s\right|\\
&\qquad \leq
\frac{\left\| \II_{\tau T}(\hat\vartheta
_{T^\delta } )^{-1} \right\|}{T^\gamma  } \left| \frac{1}{ T^{\delta -\gamma }}  \int_{T^\delta  }^{\tau T}{\dot 
  m(\hat\vartheta _{T^\delta } ,s)}\, {\rm d}\bar W_s\right|,
\end{align*}
here $\gamma $ is such that $\delta -\gamma >\frac{1}{2}$. Hence
\begin{align*}
&\Pb_{\vartheta _0}\left\{\left| \frac{\II_{\tau T}(\hat\vartheta
_{T^\delta } )^{-1}}{\tau T}  \int_{T^\delta }^{\tau T}{\dot 
  m(\hat\vartheta _{T^\delta } ,s)}\, {\rm d}\bar W_s\right|>\frac{\nu
}{4}\right\}\\
&\qquad \leq \Pb_{\vartheta _0}\left\{ \frac{1}{ T^{\delta -\gamma }} \left| \int_{T^\delta }^{\tau T}{\dot  
  m(\hat\vartheta _{T^\delta } ,s)}\, {\rm d}\bar W_s\right|>\frac{\sqrt{\nu}
}{2}\right\}\\
&\qquad\qquad  +\Pb_{\vartheta _0}\left\{ \frac{\left\|\II_{\tau T}(\hat\vartheta
_{T^\delta } )^{-1}\right\|}{ T^\gamma }  >\frac{\sqrt{\nu}
}{2}\right\}\longrightarrow 0,
\end{align*}
as $T\rightarrow \infty $, because
\begin{align*}
&\Pb_{\vartheta _0}\left\{ \frac{1}{ T^{\delta -\gamma }} \left| \int_{T^\delta }^{\tau T}{\dot  
  m(\hat\vartheta _{T^\delta } ,s)}\, {\rm d}\bar W_s\right|>\frac{\sqrt{\nu}
}{2}\right\}\\
&\qquad\qquad  \leq  \frac{1}{\nu T^{2\delta -2\gamma } }\Ex_{\vartheta _0}\int_{T^\delta }^{ T}\left|{\dot  
  m(\hat\vartheta _{T^\delta } ,s)}\right|^2\, {\rm d}s\leq \frac{C}{\nu T^{2\delta -2\gamma -1} }\rightarrow 0.
\end{align*}
Recall that $ 2\delta -2\gamma -1>0$. Further
\begin{align*}
&\Pb_{\vartheta _0}\left\{\left| \frac{\II_{\tau T}(\hat\vartheta
_{T^\delta } )^{-1}}{\tau T}  \int_{T^\delta }^{\tau T}{\dot 
  m(\hat\vartheta _{T^\delta } ,s)}\,\left[ m(\vartheta _0
  ,s)- m(\hat\vartheta _{T^\delta }
  ,s)\right]{\rm d}s\right|>\frac{\nu }{4}\right\}\\
&\quad \leq \Pb_{\vartheta _0}\left\{ \frac{1}{\tau  T^{1-\gamma }} \left| \int_{T^\delta }^{\tau T}\dot 
  m(\hat\vartheta _{T^\delta } ,s)\,\int_{0}^{1}\dot  m(\vartheta _v,s)^*
 {\rm d}v{\rm d}s\left(\hat\vartheta _{T^\delta }-\vartheta
 _0\right)\right|>\frac{\sqrt{\nu} }{2}\right\}\\
&\qquad +\Pb_{\vartheta _0}\left\{\frac{\left\|\II_{\tau T}(\hat\vartheta_{T^\delta } )^{-1}\right\|}{ T^\gamma } >\frac{\sqrt{\nu} }{2}\right\}
\longrightarrow 0
\end{align*}
as $T\rightarrow \infty $, because $\hat\vartheta _{T^\delta }-\vartheta _0=O\left(T^{-\delta /2}\right) $  and other terms are
bounded in probability. Here $\vartheta
_v=\vartheta _0+  v (\hat\vartheta _{T^\delta }-\vartheta _0 )$. 

To prove \eqref{36} we write 
\begin{align*}
&\sqrt{\tau T}\left( \vartheta
_T^\star\left(\tau \right)-\vartheta _0 \right)=\sqrt{\tau T}\left(\hat \vartheta
_{T^\delta} -\vartheta _0 \right)+\frac{\II_{\tau T}(\hat\vartheta
_{T^\delta } )^{-1}}{\sqrt{\tau T}}  \int_{T^\delta }^{\tau T}{\dot 
  m(\hat\vartheta _{T^\delta } ,s)}\, {\rm d}\bar W_s\\
&\qquad +\frac{\II_{\tau T}(\hat\vartheta
_{T^\delta } )^{-1}}{\sqrt{\tau T}}  \int_{T^\delta }^{\tau T}{\dot 
  m(\hat\vartheta _{T^\delta } ,s)}\,\left[ m(\vartheta _{0 }
  ,s)- m(\hat\vartheta _{T^\delta }
  ,s)\right]{\rm d}s.
\end{align*}
We have the estimate 
\begin{align*}
&\Ex_{\vartheta _0}\left|\frac{1}{\sqrt{\tau T}}  \int_{T^\delta }^{\tau T}\left[\dot 
  m(\hat\vartheta _{T^\delta } ,s)-\dot 
  m(\vartheta _{0} ,s)\right]\, {\rm d}\bar W_s\right|^2\\
&\qquad \leq \frac{1}{{\tau
    T}}  \int_{T^\delta }^{\tau T}\Ex_{\vartheta _0} 
\left|  \dot 
  m(\hat\vartheta _{T^\delta } ,s)-\dot 
  m(\vartheta _{0} ,s) \right|^2{\rm d}s\longrightarrow 0
\end{align*}
as $T\rightarrow \infty $, and by the central limit theorem the convergence in distribution
\begin{align*}
\frac{1}{\sqrt{\tau T}}  \int_{T^\delta }^{\tau T}\dot 
  m(\vartheta _{0} ,s)\, {\rm d}\bar W_s\Longrightarrow {\cal N}\left(0,\II(\vartheta
_{0 } ) \right). 
\end{align*}
Further, let us denote $\hat v_{T^\delta }=\sqrt{\tau T}\left(\hat \vartheta
_{T^\delta} -\vartheta _0 \right)$, then we can write 
\begin{align*}
&\hat v_{T^\delta }+\frac{\II_{\tau T}(\hat\vartheta
_{T^\delta } )^{-1}}{\sqrt{\tau T}}  \int_{T^\delta }^{\tau T}\dot 
  m(\hat\vartheta _{T^\delta } ,s)\left[{
  m(\vartheta _{0 } ,s})-  
  m(\hat\vartheta _{T^\delta } ,s) \right]\, {\rm d}s\\
&\quad =\II_{\tau T}(\hat\vartheta
_{T^\delta } )^{-1}\left(  \II_{\tau T}(\hat\vartheta
_{T^\delta } )-\frac{1}{\tau T}\int_{0 }^{1}\int_{T^\delta }^{\tau T}\dot 
  m(\hat\vartheta _{T^\delta } ,s)\dot 
  m(\vartheta _{r } ,s)^*{\rm d}r\;{\rm d}s\right)\hat v_{T^\delta },
\end{align*}
here $\vartheta _r= \hat\vartheta _{T^\delta }+r \left(\hat\vartheta
_{T^\delta }-\vartheta_0\right)  $. The presentation
\begin{align*}
\dot   m(\vartheta _{r } ,s)=\dot   m(\hat\vartheta_{T^\delta } ,s)+\int_{0
}^{1}\ddot   m(\vartheta_{q } ,s){\rm d}q \left(\hat\vartheta
_{T^\delta }-\vartheta_0\right) 
\end{align*}
and the equality 
\begin{align*}
\II_{\tau T}(\hat\vartheta_{T^\delta } )=\frac{1}{\tau T}\int_{T^\delta }^{\tau T}\dot 
  m(\hat\vartheta _{T^\delta } ,s)\dot 
  m(\hat\vartheta _{T^\delta } ,s)^*{\rm d}s
\end{align*}
    allows us to write
\begin{align*}
&\hat v_{T^\delta }+\frac{\II_{\tau T}(\hat\vartheta
_{T^\delta } )^{-1}}{\sqrt{\tau T}}  \int_{T^\delta }^{\tau T}\dot 
  m(\hat\vartheta _{T^\delta } ,s)\left[{
  m(\vartheta _{0 } ,s})-  
  m(\hat\vartheta _{T^\delta } ,s) \right]\, {\rm d}s\\
&\qquad \quad =\sqrt{\tau
    T}\left|\hat\vartheta_{T^\delta } -\vartheta _0 \right|^2O\left( 1\right)
=T^{\frac{1}{2}-\delta }O\left( 1\right)\longrightarrow 0,
\end{align*}
as $T\rightarrow \infty $,

Let us verify that the Fisher information matrix is non degenerate. It is
sufficient to show that the matrix 
\begin{align*}
\JJ\left({\vartheta _0} \right)=\begin{pmatrix}
\Ex_{\vartheta _0} \tilde y_t^2,  &\Ex_{\vartheta _0}\tilde y_t\tilde z_t \\
\Ex_{\vartheta _0} \tilde y_t\tilde z_t,  &\Ex_{\vartheta _0} \tilde z_t^2
\end{pmatrix},
\end{align*}
is non degenerated. Here $\tilde y_t, \tilde z_t $ are stationary solutions of \eqref{yt} and
\eqref{zt} respectively. If this matrix is degenerated,  then
\begin{align}
\label{c-s}
\Ex_{\vartheta _0}\tilde y_t^2\;\Ex_\vartheta\tilde z_t^2  =\left(\Ex_{\vartheta _0} \tilde y_t\tilde z_t \right)^2.
\end{align}
Recall that by Cauchy-Schwarz inequality
\begin{align*}
\left(\Ex_{\vartheta _0} \tilde y_t\tilde z_t \right)^2\leq \Ex_{\vartheta _0} \tilde
y_t^2\;\Ex_{\vartheta _0}\tilde z_t^2  
\end{align*}
with equality if and only if  $\tilde z_t=c\tilde y_t$ with some constant $c\not
=0$. Therefore in the case of equality we have $\Ex_{\vartheta _0}\left(c\tilde y_t-\tilde z_t \right)^2=0 $.

Introduce a new process $\tilde v_t=c\tilde y_t-\tilde z_t$ as a solution of the equation
\begin{align*}
{\rm d}\tilde v_t&=\left[\tilde  x_t\left(1-c\right) -1 \right] {\rm d}t-\left[a+b^2\tilde x_t\left(1-\tilde x_t\right)\right]\tilde v_t{\rm
  d}t+b\left(1-2\tilde x_t\right)\tilde v_t\;{\rm d}\bar W_t,
\end{align*}
here $\tilde   v_t$ and $\tilde   x_t$  are stationary solutions. 

Further, following \cite{PCh09}, Section 4,  here the similar estimate was obtained, we
write this solution as
\begin{align*}
\tilde v_t&=\tilde v_0e^{-at}+\int_{0}^{t}e^{-a\left(t-s\right)} \left[\tilde
  x_s\left(1-c\right) -1 \right]\; {\rm d}s -b^2\int_{0}^{t}e^{-a\left(t-s\right)} \tilde
x_s\left(1-\tilde x_s\right)\tilde v_s\;{\rm d}s\\
&\quad + b\int_{0}^{t}e^{-a\left(t-s\right)} \left(1-2\tilde x_s\right)\tilde
v_s\;{\rm d}\bar W_s.
\end{align*}
Hence
\begin{align*}
&\Ex_{\vartheta _0}\left( \int_{0}^{t}e^{-a\left(t-s\right)} \left[\tilde
  x_s\left(1-c\right) -1 \right] {\rm d}s  \right)^2\leq
4\left(1+e^{-2at}\right)\Ex_{\vartheta _0}\tilde v_t^2 \\
&\qquad  +\frac{4b^4}{a}\int_{0}^{t}e^{-a\left(t-s\right)} \frac{1}{16}
\Ex_{\vartheta _0} \tilde v_s^2 \;{\rm d}s
+4b^2\int_{0}^{t}e^{-2a\left(t-s\right)}\Ex_{\vartheta _0}\tilde v_s^2 \;{\rm d}s
\leq C \Ex_{\vartheta _0} \tilde v_t^2
\end{align*}
with some constant $C>0$ which does not depend on $t$. Recall that
$\Ex_{\vartheta _0} \tilde v_t^2 $ does not depend on $t$ too because $\tilde
v_t $ is stationary solution. Therefore if we  show that for all $c$
\begin{align*}
\lim_{t\rightarrow \infty }\Ex_{\vartheta _0}\left( \int_{0}^{t}e^{-a\left(t-s\right)} \left[\tilde
  x_s\left(1-c\right) -1 \right] {\rm d}s  \right)^2>0,
\end{align*}
then the matrix $\JJ\left(\vartheta _0\right)$ is non degenerate. 
The random process
\begin{align*}
\zeta _t=\int_{0}^{t}e^{-a\left(t-s\right)} \left[\tilde
  x_s\left(1-c\right) -1 \right] {\rm d}s  
\end{align*}
is the solution of the equation
\begin{align*}
\frac{{\rm d}\zeta _t}{{\rm d}t}=-a\zeta _t+\tilde
  x_t\left(1-c\right) -1 ,\qquad \zeta _0=0. 
\end{align*}
The elementary calculations show that for all $\vartheta _0$ and $c$
 \begin{align*}
\lim_{t\rightarrow \infty }\Ex_{\vartheta _0}\zeta _t^2=\frac{\Ex_{\vartheta _0} \left[\tilde
  \pi _0\left(1-c\right) -1\right]^2 }{a^2}>0,
\end{align*}
here $\tilde \pi _0 $ is the stationary distribution. 

Therefore the Fisher information matrix  $\II\left(\vartheta _0\right)$ is non
degenerate for all $\vartheta_0\in \Theta  $.

Note that the limit covariance matrix of the one-step MLE-process by the
Theorem \ref{T2} coincides
with the covariance of the asymptotically efficient MLE \cite{PCh09}, therefore  $\vartheta
_{T}^\star \left(\tau \right)$ is asymptotically efficient too.

\section{Discussions}

The learning interval in one-step section is $[0,T^\delta ]$, where
$\delta \in (\frac{1}{2}, 1)$, i.e., it is negligeable with respect to the
whole observations time $T$. It can be done even shorter, if we use two-step
MLE-process approach, as it was proposed in \cite{Kut15}. It corresponds to
the learning interval $[0,T^\delta )$ with $\delta \in (\frac{1}{4},
  \frac{1}{2}]$. The procedure is the follows. First we obtain the preliminary
estimator $\hat\vartheta _{T^\delta }$ as before. Then we introduce the second
preliminary estimator
\begin{align*}
\vartheta _{t,T}^\star=\hat\vartheta _{T^\delta
}+t^{-1}\II_t(\hat\vartheta _{T^\delta })^{-1} \Delta
_t(\hat\vartheta _{T^\delta },X^t ) ,\qquad t\in \left[T^\delta ,T\right]
\end{align*}
and then we define two-step MLE-process
\begin{align*}
\vartheta _{t,T}^{\star\star}=\vartheta _{t,T}^\star+t^{-1}\II_t(\hat\vartheta _{T^\delta })^{-1} \Delta
_t(\hat\vartheta _{T^\delta },\vartheta _{t,T}^\star,X^t ) ,\qquad t\in \left[T^\delta ,T\right],
\end{align*}
here
\begin{align*}
\Delta
_t\left(\vartheta _1,\vartheta _2,X^t \right)=\frac{1}{\sqrt{t}}
\int_{T^\delta }^{t} \dot m\left(\vartheta _1,s\right)\left[{\rm
    d}X_s-m\left(\vartheta _2,s\right) {\rm d}s\right]   ,\qquad t\in \left[T^\delta ,T\right].
\end{align*}
It can be shown that for all $\tau \in (0,1]$ and $t=\tau T$ we have the
  asymptotic normality of the estimator $\vartheta _{T}^{\star\star}\left(\tau
  \right)=\vartheta _{\tau T,T}^{\star\star} $: 
\begin{align*}
\sqrt{\tau T}\left( \vartheta _{T}^{\star\star}\left(\tau \right)-\vartheta_0
\right)\Longrightarrow {\cal N}\left(0, \II\left(\vartheta_0 \right)^{-1}\right).
\end{align*}
See the details in \cite{Kut15}.
\bigskip

Note that it can be shown that the one-step MLE-process converges in
distribution to the limit Brownian motion. Let us  introduce the random process
\begin{align*}
\eta _T\left(\tau \right)=\tau \sqrt{T}\II(\vartheta_0 )^{-1/2}\left(\vartheta
_{ T}^\star\left(\tau \right)-\vartheta _0\right) ,\qquad \tau _\delta \leq \tau \leq 1,
\end{align*}
here  $\tau _\delta =T^{\delta -1}\rightarrow 0$. 
More detailed analysis shows that the random process $\eta _T\left(\tau
\right), \tau _*\leq \tau \leq 1 $ converges to two-dimensional standard
Wiener process $W\left(\tau \right), \tau _*\leq \tau \leq 1$ with any $\tau
_*\in (0,1]$. For the details see the proof of such convergence in similar
  problem in \cite{Kut15}.

\bigskip

{\bf Acknowledgment.} This work was done under partial financial support
(second author) of
the grant of  RSF number 14-49-00079.


\begin{thebibliography}{99}
\bibitem {BRR98} Bickel, P.J., Ritov, Y. and Ryd\'en, T. (1998) Asymptotic
  normality of the maximum likelihood estimator for general hidden Markov
  models. {\it Ann. Statist.}, 26, 4, 1614-1635.
\bibitem {CMT05}  Capp\'e, O., Moulines, E. and Ryd\'en, T. (2005) {\it Inference in Hidden
  Markov Models}. Springer, N.Y.
\bibitem {PCh09} Chigansky, P. (2009) Maximum likelihood estimation for hidden
  Markov models in continuous time. {\it Statist. Inference Stoch. Processes},
  12, 2, 139-163.
\bibitem {DZ86}  Dembo, A. and Zeitouni, O. (1986) Parameter estimation for
  partially observed continuous time processes via the EM algorithm. {\it
    Stoch. Proces. Applic.} 23, 91-113.
\bibitem {EAM95} Elliott, R.J., Aggoun, L. and Moor, J.B. (1995) {\it Hidden
  Markov Models.} Springer,  N.Y.
\bibitem {IY09} Iacus, S. and Yoshida, N. (2009) Estimation for discretly
  observed telegraph process. {\it Theory Probab. Math. Statist.}, 78, 37-47.
\bibitem {IH81} Ibragimov I.A. and Khasminskii R. (1981) {\it Statistical
Estimation - Asymptotic Theory.} {Springer-Verlag}, New York.
\bibitem{KU15} Kamatani, K. and Uchida, M. (2015) {Hybrid multi-step
  estimators for stochastic 
  differential equations based on sampled data.}
   {\it  Statist. Inference  Stoch. Processes.} 18, 2, 177-204.
\bibitem{K66}  Khasminskii, R.Z. (1966)   On stochastic processes defined by
differential equations with small parameter. {\sl Theory Probab. Appl.,}
{\bf 11}, 240--259.
\bibitem {Kut84} Kutoyants, Yu.A. (1984) {\it Parameter Estimation for
  Stochastic Processes.} Heldermann, Berlin. 
\bibitem{Kut94}  Kutoyants, Yu.A. (1994) {\it Identification of
        Dynamical Systems with Small Noise.} Kluwer,  Dordrecht.
\bibitem {Kut14} Kutoyants, Yu.A. (2014) { On approximation of the backward
  stochastic differential equation. Small noise, large samples and high
  frequency cases.}  {\it Proceedings of the Steklov Institute of
    Mathematics.}, 287, 133-154.
\bibitem {Kut15} Kutoyants, Yu.A. (2015) On multi-step MLE-processes for
  ergodic diffusion. submitted.
\bibitem {LS01} Liptser, R.S. and Shiryaev, A.N. (2001) {\it Statistics of
  Random Processes}, 2-nd ed., vol. 1, Springer, N.Y.
\bibitem {V87} Veretennikov, A. Yu. (1987) Bounds for the mixing rate in the
  theory of stochastic equations. {\it Theory Probab. Appl.}, 32, 2, 273-281.
\bibitem {W65} Wonham, W. M. (1965) Some applications of stochastic differential
  equations to optimal non-linear filtering. {\it SIAM J. Contr.}, Ser. A, 2,
  347-369. 
\bibitem {ZD88} Zeitouni, O.  and Dembo, A.  (1988) Exact filters for the
  estimation of the number of transitions of finite-state continuous time
  Marcov processes. {\it IEEE IT}, 34, 4, 890-893
\end{thebibliography}
\end{document}